\documentclass[12pt]{article} 
\newcommand{\be}{\begin{equation}}
\newcommand{\ee}{\end{equation}}
\newcommand{\ba}{\begin{array}}
\newcommand{\ea}{\end{array}}

\newcommand{\bea}{\begin{eqnarray*}}
\newcommand{\eea}{\end{eqnarray*}}
\newcommand{\bean}{\begin{eqnarray}}
\newcommand{\eean}{\end{eqnarray}}
\newcommand{\proof}{\vspace{1ex}\noindent{\em Proof}. \ }
\def\ds{\displaystyle}
\def\nm{\noalign{\medskip}}

\newcommand{\grad}{{\rm \; grad\;}}

\newtheorem{lemma}{Lemma}[section]

\newtheorem{theorem}{Theorem}[section]
\newtheorem{corollary}{Corollary}[section]
\newtheorem{proposition}{Proposition}[section]
\usepackage{amsmath, amssymb}
\usepackage{latexsym}
\newcommand{\R}{\mathbb{R}}

\begin{document}
\title{On a hyperbolic coefficient inverse problem via partial dynamic boundary measurements}
\author{Christian Daveau \thanks{
D\'epartement de Math\'ematiques, CNRS AGM UMR 8088, Universit\'e de
Cergy-Pontoise, 95302 Cergy-Pontoise Cedex, France (Email:
christian.daveau@math.u-cergy.fr).} \and Abdessatar Khelifi
\thanks{ D\'epartement de Math\'ematiques, Universit\'e des Sciences
de Carthage, Bizerte, Tunisia. (Email:
abdessatar.khelifi@fsb.rnu.tn).}}

\maketitle
\maketitle \abstract{ This paper is devoted to the identification of
the unknown smooth coefficient c entering the hyperbolic equation
$c(x)\partial_{t}^{2}u - \Delta u = 0$ in a bounded smooth domain in
$\R^{d}$ from partial (on part of the boundary) dynamic boundary
measurements. In this paper we prove  that the knowledge of the
partial Cauchy data for this class of hyperbolic PDE on any open
subset $\Gamma$ of the boundary determines explicitly the
coefficient $c$ provided that $c$ is known outside a bounded domain.
Then, through construction  of appropriate test functions by a
geometrical control method, we derive a formula for calculating the
coefficient $c$ from the knowledge of the difference between the
local Dirichlet to Neumann maps.}



\section{Introduction}
In this paper, we present a differently  method for multidimensional
Coefficient Inverse Problems (CIPs) for a class of hyperbolic
Partial Differential Equations (PDEs). In the literature, the reader
can find many key investigations in this kind of inverse problems,
see, e.g.
\cite{ammari,b-kalibanov1,belishev1,belishev2,Kab,klibanov1,nachman,Ramm1,RRa,Y,Y95}
and references cited there. L. Beilina and M.V. Klibanov have deeply
studied this important problem in various recently works
\cite{b-kalibanov1, b-kalibanov2}. In \cite{b-kalibanov1}, the
authors have introduced a new globally convergent numerical method
to solve a coefficient inverse problem associated to a hyperbolic
PDE. The development of globally convergent numerical methods for
multidimensional CIPs has started, as a first generation, from the
developments found in
\cite{klibanov-timinov2007,klibanov-timinov2004,klibanov-xin}. Else,
A. G. Ramm and Rakesh have developed a general method for proving
uniqueness theorems for multidimensional inverse problems. For the
two dimensional case, Nachman \cite{nachman} proved an uniqueness
result for CIPs for some elliptic equation. Moreover, we find the
works of L. P$\ddot{a}$iv$\ddot{a}$rinta and V. Serov
\cite{serov-p1,serov-p2} about the same issue, but for elliptic
equations.
\\ In
other manner, the author Y. Chen has treated in \cite{chen} the
Fourier transform of the hyperbolic equation similar to ours with
the unknown coefficient $c(x)$. Unlike this, we derive, using as
weights particular background solutions constructed by a geometrical
control method, asymptotic formulas in terms of the partial dynamic
boundary measurements (Dirichlet-to-Neumann map) that are caused by
the small perturbations. These asymptotic formulae yield the inverse
Fourier transform of unknown coefficient.\\

The ultimate objective of the work described in this paper is to
determine, effectively, the unknown smooth coefficient c entering
 a class of hyperbolic equation in a bounded smooth domain in $\R^{d}$ from partial (on part of the
boundary) dynamic boundary measurements. The main difficulty appears
in boundary measurements, is that the formulation of our boundary
value problem involves unknown boundary values. This problem is well
known in the study of the classical elliptic equations, where the
characterization of the unknown Neumann boundary value in terms of
the given Dirichlet datum is known as the Dirichlet-to-Neumann map.
But, the problem of determining the unknown boundary values also
occurs in the study of hyperbolic equations formulated in a bounded
domain.

 As our main result we develop, using
as weights particular background solutions constructed by a
geometrical control method, asymptotic formulas for appropriate
averaging of the partial dynamic boundary measurements that are
caused by the small perturbations of coefficient according to a
parameter $\alpha$. Assume that the coefficient is known outside a
bounded domain $\Omega$, and suppose that we know explicitly the
value of $\displaystyle\lim_{\alpha\to 0^+}c(x)$ for $x\in\Omega$.
Then, the developed asymptotic formulae yield the inverse Fourier
transform of
the unknown part of this coefficient.\\

 In the subject of small volume perturbations from a known background material associated to the full time-dependent Maxwell's
 equations, we have
derived asymptotic formulas to identify their locations and certain
properties of their shapes from dynamic boundary measurements
\cite{dav-khe-sus}. The present paper represents a different
investigation of this line of work.

As closely related stationary identification problems we refer the
reader to \cite{CMV,FV,nachman,SU} and references cited there.

\section{Problem formulation}

Let $\Omega$ be a bounded, smooth subdomain of $\R^d$ with $d\leq3$,
(the assumption $d\leq3$ is necessary in order to obtain the
appropriate regularity for the solution using classical Sobolev
embedding, see Brezis \cite{brezis}). For simplicity we take
$\partial \Omega$ to be ${\cal C}^{\infty}$, but this condition
could be considerably weakened. Let $n=n(x)$ denote the outward unit
normal vector to $\Omega$ at a point on $\partial \Omega$. Let
$T>0$, $x_0\in \R^d \Omega$ and let $\Omega^\prime$ be a smooth
subdomain of $\Omega$. We denote by $\Gamma \subset \subset
\partial \Omega$ as a measurable smooth open part of the boundary
$\partial \Omega$.\\
Throughout this paper we shall use quite standard $L^2-$ based
Sobolev spaces to measure regularity.

As the forward problem, we consider the Cauchy problem for a
hyperbolic PDE
\begin{eqnarray}\label{wave1}
  c(x)v_{tt}-\Delta v=0& \mbox{in }\R^d\times (0,T) \\
  v(x,0)=0, & v_t(x,0)=\delta (x-x_0)+\chi(\Omega)\psi,
\end{eqnarray}
where $\chi(\Omega)$ is the characteristic function of $\Omega$ and
$\psi \in {\cal
C}^{\infty}(\R^d)$ that $\psi(x)\neq0,$ $\forall x\in \overline{\Omega}$.\\

Equation (\ref{wave1}) governs a wide range of applications,
including e.g., propagation of acoustic and electromagnetic waves.\\
We assume that the coefficient $c(x)$ of equation (\ref{wave1}) is
such that
\begin{equation}\label{coeff}
  c(x)=\begin{cases}
  c_\alpha(x)=c_{0}(x)+\alpha c_{1}(x) & \mbox{for }x\in\Omega, \\
  c_2(x)=const.>0 & \mbox{for }x\in\R^d\backslash\Omega;
  \end{cases}
\end{equation}
where $c_i(x)\in {\cal C}^2(\overline{\Omega})$ for $i=0,1$ with
\begin{equation}\label{c0-c1}
c_1 \equiv 0 {\rm \;in\;} \Omega \setminus
\overline{\Omega^\prime},\quad \mbox{and }M:=sup
\{c_1(x);x\in\Omega^\prime \},\end{equation} where $\Omega^\prime$
is a smooth subdomain of $\Omega$ and $M$ is a positive constant. We
also assume that $\alpha
>0$, the order of magnitude of the small perturbations of coefficient, is sufficiently small
 that \begin{equation} \label{ineq1} |c_\alpha(x)|
\geq c_* > 0, x \in \overline{\Omega}, \end{equation} where $c_*$ is
a positive constant.\\

Suppose that the positive number $c_2$ is given. In this paper we
assume that the function $c(x)$ is unknown in the domain $\Omega$.
Our purpose is the determination of $c(x)$ for $x\in\Omega$,
assuming that the following function $g(x, t)$ is known for the
single source position $x_0 \in\R^d\backslash\overline{\Omega}$.
Therefore, as done for the Dirichlet boundary conditions
in~\cite{b-kalibanov2}, we set the Neumann boundary conditions:
\begin{equation}\label{bound-g}
\frac{\partial v}{\partial n}|_{\partial\Omega\times(0,T)}=g(x,t).
\end{equation}
The knowledge of $c(x)$ outside of $\Omega$ ($c(x)=c_2$ in
$\R^d\backslash \Omega$), and the boundary function $g(x,t)$ allow
us to determine uniquely the function $v(x,t)$ for $x\in
\R^d\backslash \Omega$ as solution of the boundary value problem for
equations (\ref{wave1})-(2) with initial conditions in (2) and with
the boundary conditions (\ref{bound-g}). Therefore, one can uniquely
determine the function $f(x,t)=v|_{\partial\Omega\times(0,T)}$.

Then, we can now consider an initial boundary value problem only in
the domain $\Omega\times(0,T)$. Thus, the function $v$ satisfying
(\ref{wave1})-(2) is en particular solution of the following initial
boundary value problem
 \begin{equation} \label{walpha} \left\{
\begin{array}{l}
\ds (c_\alpha\partial_t^2 -\Delta )u_\alpha= 0 \quad {\rm in}\;
\Omega \times (0, T),\\ \nm u_\alpha |_{t=0} = \varphi,
\partial_t u_\alpha |_{t=0} = \psi \quad {\rm in}\; \Omega,\\ \nm
u_\alpha |_{\partial \Omega \times (0, T)} = f.
\end{array}
\right. \end{equation}

Define $u$ to be the solution of the hyperbolic equation in the
homogeneous situation ($\alpha=0$). Thus, $u$ satisfies
\begin{equation} \label{wo} \left\{
\begin{array}{l}
\ds (c_0\partial_t^2 - \Delta) u = 0 \quad {\rm in}\;
\Omega \times (0, T),\\
\nm u |_{t=0} = \varphi, \partial_t u |_{t=0} = \psi \quad {\rm
in}\;
\Omega,\\
\nm u |_{\partial \Omega \times (0, T)} = f.
\end{array}
\right. \end{equation}
 Here  $\varphi\in {\cal
C}^{\infty}(\overline{\Omega})$
 and $f \in {\cal C}^{\infty}(0, T; {\cal C}^{\infty}(\partial
 \Omega))$ are subject to the compatibility conditions
$$ \partial_{t}^{2l} f |_{t = 0} = (\Delta^l \varphi)
|_{\partial \Omega} {\rm \; and \;} \partial_t^{2l + 1}  f |_{t =0}
= (\Delta^l \psi ) |_{\partial \Omega}, \quad l=1, 2, \ldots $$
which give that (\ref{wo}) has a unique solution in ${\cal C}^\infty
([0, T] \times \overline{\Omega})$, see \cite{E}. It is also
well-known that (\ref{walpha}) has a unique weak solution $u_\alpha
\in {\cal C}^0(0, T; H^1(\Omega)) \cap {\cal C}^1(0, T;
L^2(\Omega))$, see \cite{L}, \cite{E}. Indeed, from \cite{L} we have
that $\ds \frac{\partial u_\alpha}{\partial n} |_{\partial \Omega} $
belongs to $L^2(0, T; L^2(\partial \Omega))$.\\

Now, we define $\Gamma_c:=\partial \Omega \setminus
\overline{\Gamma}$, and we introduce the trace space
$$\widetilde{H}^{\frac{1}{2}}(\Gamma) = \Bigr\{ v \in
{H}^{\frac{1}{2}}(\partial \Omega\times(0,T)), v \equiv 0 \mbox{ on
} \Gamma_c \times(0,T)\Bigr\}.$$  It is known that the dual of
$\widetilde{H}^{\frac{1}{2}}(\Gamma)$ is
${H}^{-\frac{1}{2}}(\Gamma)$.

To introduce the local Dirichlet to Neumann map associated to our
problem, we firstly define the function $\tilde{f}=\chi(\Gamma)
f\quad \mbox{for } (x,t)\in\partial\Omega\times(0,T)$, where
$\chi(\Gamma)$ is the characteristic function of $\Gamma$. Then, we
have
\begin{equation}\label{f-tilde}
\tilde{f}=f|_{\Gamma\times(0,T)} \mbox{ and }\tilde{f}\in
\widetilde{H}^{\frac{1}{2}}(\Gamma).
\end{equation}
Therefore, we define the local Dirichlet to Neumann map associated
to coefficient $c_\alpha$ by
:\[\displaystyle\Lambda_{\alpha}(\tilde{f}) = \frac{\partial
u_\alpha}{\partial n}|_{\Gamma} \mbox{ for }\tilde{f} \in
\widetilde{H}^{\frac{1}{2}}(\Gamma),\] where $u_\alpha$ is the
solution of (\ref{walpha}). Let $u$ denote the solution to the
hyperbolic equation (\ref{wo}) with the Dirichlet boundary condition
$u = f$  on $\partial\Omega\times(0,T)$. Then,  the local Dirichlet
to Neumann map associated to $c_0$ is
$\displaystyle\Lambda_{0}(\tilde{f}) = \frac{\partial u}{\partial
n}|_{\Gamma}$ for $\tilde{f} \in
\widetilde{H}^{\frac{1}{2}}(\Gamma).$\\

Our problem can be stated as follows:\\

 \textbf{Inverse problem.}
Suppose that the smooth coefficient $c(x)$ satisfies (3)-(4)-(5),
where the positive number $c_2$ is given. Assume that the function
$c(x)$ is unknown in the domain $\Omega$ and $\tilde{f}$ is given by
(\ref{f-tilde}). Is it possible to determine the coefficient
$c_{\alpha}(x)$ from the knowledge of the difference between the
local Dirichlet to Neumann maps $\ds
\Lambda_{\alpha} - \Lambda_0$ on $\Gamma$, if
we know explicitly the value of
$\displaystyle\lim_{\alpha\to 0^{+}c_\alpha}(x)$ for $x\in\Omega$ ? \\

To give a positive answer, we will develop an asymptotic expansions
of an "appropriate averaging" of $ \ds \frac{\partial
u_\alpha}{\partial n}$ on $ \Gamma \times (0, T)$, using particular
background solutions as weights. These particular solutions are
constructed by a control method as it has been done in the original
work \cite{Y} (see also \cite{belishev-k}, \cite{BY}, \cite{PY},
\cite{PY1} and \cite{Y95}). It has been known for some time that the
full knowledge of the (hyperbolic) Dirichlet to Neumann map
($u_\alpha |_{\partial \Omega \times (0, T)} \mapsto \frac{\partial
u_\alpha}{\partial n} |_{\partial \Omega \times (0, T)}$) uniquely
determines conductivity, see \cite{RS}, \cite{S}. Our identification
procedure can be regarded as an important attempt to generalize the
results of \cite{RS} and \cite{S}
 in the case of partial knowledge (i.e., on only part of the boundary) of the Dirichlet
 to Neumann map to determine the coefficient of the hyperbolic equation considered
 above. The question of uniqueness of this inverse problem can be addressed
positively via the method of Carleman estimates, see, e.g., \cite{klibanov1,klibanov-timinov2004}.\\

\section{The Identification Procedure}
Before describing our identification procedure, let us introduce the
 following cutoff function $\beta(x) \in
{\cal C}^{\infty}_0(\Omega)$ such that
$\beta \equiv 1$ on $\Omega^{\prime}$ and let $\eta \in \R^d$.\\
We will take in what follows $\varphi(x) = e^{i \eta \cdot x},
\psi(x) =  - i |\eta | e^{i \eta \cdot x}, {\rm \; and \;} f(x, t) =
e^{i \eta \cdot x  -  i | \eta | t}$ and assume that we are in
possession of the boundary measurements of
$$\ds \frac{\partial u_\alpha}{\partial n} \quad {\rm on\;} \Gamma
\times (0, T).$$
 This particular choice of data $\varphi, \psi,$ and $f$
implies that the background solution $u$ of the wave equation
(\ref{wo}) in the homogeneous background medium can be given
explicitly.\\

Suppose now that $T$ and the part $\Gamma$ of the boundary $\partial
\Omega$ are such that they  geometrically control $\Omega$ which
roughly means that every geometrical optic ray, starting at any
point $x \in \Omega$ at time $t=0$ hits $\Gamma$ before time $T$ at
a non diffractive point, see \cite{BLR}. It follows from \cite{Y_98}
(see also \cite{AW}) that there exists (a unique) $g_\eta \in
H^1_0(0, T; TL^2(\Gamma))$ (constructed by the Hilbert Uniqueness
Method) such that the unique weak solution $w_\eta$ to the wave
equation \begin{equation} \label{weta} \left\{
\begin{array}{l}
\ds (c_0\partial_t^2 -  \Delta) w_\eta  = 0 \quad {\rm in}\;
\Omega \times (0, T),\\
\nm
w_\eta  |_{t=0} = \beta(x) e^{i \eta \cdot x} \in H^1_0(\Omega), \\
\nm
\partial_t w_\eta
 |_{t=0} = 0 \quad {\rm in}\;
\Omega,\\
\nm
w_\eta |_{\Gamma  \times (0, T)} = g_\eta,\\
\nm w_\eta |_{\partial \Omega \setminus \overline{\Gamma}
  \times (0, T)} = 0,
\end{array}
\right. \end{equation} satisfies $w_\eta(T) = \partial_t w_\eta (T)
= 0$.

Let $\theta_\eta \in H^1(0, T; L^2(\Gamma))$ denote the unique
solution of the Volterra equation of second kind
\begin{equation}\label{eq4m} \left\{\begin{array}{l} \ds
\partial_t \theta_\eta (x, t) + \int_t^T e^{- i | \eta | (s -t)}
( \theta_\eta (x, s) - i | \eta | \partial_t \theta_\eta (x, s)) \;
ds = g_\eta(x, t) \quad {\rm for\;} x \in \Gamma, t \in (0, T),
\\ \theta_\eta(x, 0) = 0 \quad {\rm for\;} x \in \Gamma.
\end{array}
\right. \end{equation} We can refer to the work of Yamamoto in
\cite{Y95} who conceived the idea of using such Volterra equation to
apply the geometrical control for solving
inverse source problems.\\

The existence and uniqueness of this $\theta_\eta$ in $H^1(0, T;
L^2(\Gamma))$ for any $\eta \in \R^d$ can be established using the
resolvent kernel. However, observing from differentiation of
(\ref{eq4m}) with respect to $t$ that $\theta_\eta$ is the unique
solution of the ODE: \begin{equation} \label{eq4p} \left\{
\begin{array}{l}
\ds
\partial_t^2 \theta_\eta - \theta_\eta = e^{i |\eta| t}
\partial_t ( e^{-i |\eta| t} g_\eta) \quad {\rm for\;} x \in
\Gamma, t \in (0, T), \\ \theta_\eta(x, 0) = 0,
\partial_t \theta_\eta(x, T) = 0  \quad {\rm for\;} x \in \Gamma,
\end{array}
\right. \end{equation} the function $\theta_\eta$ may be find (in
practice) explicitly with variation of parameters and it also
immediately follows from this observation that $\theta_\eta$ belongs
to
$H^2(0, T; L^2(\Gamma))$.\\
We introduce $v_{ \eta}$ as the unique weak solution (obtained by
transposition) in $ {\cal C}^0(0, T; L^2(\Omega)) \cap {\cal C}^1(0,
T; H^{-1}(\Omega))$ to the wave equation
\begin{equation}\label{v-eta}
\left\{\begin{array}{l} \ds (c_0\partial_t^2 - \Delta) v_{\eta}  = 0
\quad {\rm in}\;
\Omega \times (0, T),\\
\nm v_{\eta}  |_{t=0} = 0 \quad {\rm in}\;
\Omega,\\
\nm \ds
 \partial_t v_{\alpha,\eta}
 |_{t=0} = i \nabla\cdot (\eta c_1(x) e^{i \eta \cdot x}) \in L^2(\Omega),\\
\nm v_{\eta}  |_{\partial \Omega  \times (0, T)} = 0.
\end{array}
\right.
\end{equation}
Then, the following holds.

\begin{proposition}\label{p4.1}
Suppose that $\Gamma$ and $T$ geometrically control $\Omega$. For
any $\eta \in \R^d$ we have \begin{equation}\label{rel-prop3.1} \ds
\int_0^T \int_\Gamma g_\eta \Lambda_{0}( v_{
\eta})~d\sigma(x)dt=|\eta|^2 \int_{\Omega^\prime}
 c_1(x) e^{2i \eta \cdot x} \; dx. \end{equation}  Here $d\sigma(x)$ means an elementary
surface for $x\in \Gamma$.
\end{proposition}
\proof Let $v_{\eta}$ be the solution of (\ref{v-eta}). From
\cite{L} [Theorem 4.1, page 44] it follows that
 $\ds \Lambda_{0}( v_{
\eta})=\frac{\partial v_{\eta}}{\partial n} |_{\Gamma} \in L^{2}(0,
T; L^2(\Gamma))$. Then, multiply the equation $\ds (\partial_t^2 +
\Delta ) v_\eta  = 0$ by $w_\eta$ and integrating by parts over $(0,
T) \times \Omega$, for any $\eta \in \R^d$ we have
\[
\ds \int_0^T \int_\Omega (\partial_t^2 -  \Delta) v_{\eta}  w_\eta =
i \int_\Omega \nabla\cdot(\eta c_1(x) e^{i \eta \cdot x}) \beta(x)
e^{i \eta \cdot x} \; dx - \int_0^T \int_\Gamma g_\eta
\frac{\partial v_{ \eta}}{\partial n} = 0.
\]
Therefore \begin{equation} \label{r1} \ds   |\eta|^2
\int_{\Omega^\prime}
 c_1(x) e^{2i \eta \cdot x} \; dx
=  \int_0^T \int_\Gamma g_\eta \frac{\partial v_{ \eta}}{\partial
n}. \end{equation} since $c_1 \equiv  0$ on $\Omega \setminus
\overline{\Omega^\prime}$.  $\square$

In term of the function $v_{\eta}$ as solution of (\ref{eq4m}), we
introduce
\begin{equation}\label{tilde-v}
\ds \tilde{u}_{\alpha}(x, t) = u(x, t) + \alpha^d \int_0^t e^{- i |
\eta | s} v_{\eta}(x, t-s)\; ds, x \in \Omega, t \in (0, T).
\end{equation}
Moreover, for  $z(t)\in {\cal C}^\infty_0(]0, T[)$ and for any $v
\in L^1(0, T; L^2(\Omega))$, we define
\begin{equation}\label{hat-v}
\ds \hat{v}(x) = \int_0^T v(x, t) z(t) \; dt \in L^2(\Omega).
\end{equation}

 The following lemma is useful to proof our main result.
\begin{lemma}\label{lemm1}Consider an arbitrary function $c(x)$
satisfying condition (\ref{coeff}) and assume that conditions
(\ref{c0-c1}) and (\ref{ineq1}) hold. Let $u$, $u_\alpha$ be
solutions of (\ref{wo}) and (\ref{walpha}) respectively. Then, using
(\ref{tilde-v}) the following estimates hold:
\[
|| u_\alpha - u ||_{L^\infty(0, T; L^2(\Omega))}  \leq C \alpha,
\]
where $C$ a positive constant. And,
\begin{equation} \label{r5} || \tilde{u}_\alpha - u_\alpha
 ||_{L^\infty(0, T; L^2(\Omega))}  \leq C^{\prime} \alpha^{d+1},
\end{equation}
where $C^{\prime}$ is a positive constant.
\end{lemma}
\proof Let $y_\alpha$ be defined by
\[
\left\{ \begin{array}{l} y_\alpha \in H^1_0(\Omega),\\ \nm \Delta
y_\alpha = c_{\alpha}\partial_t (u_\alpha - u) \quad {\rm in\;}
\Omega.
\end{array}
\right.
\]
We have
\[
\ds \int_{\Omega} c_{\alpha}\partial_t^2 (u_\alpha -u) y_\alpha
 + \int_{\Omega}
 \nabla (u_\alpha -u) \cdot \nabla y_\alpha  = \alpha
\int_\Omega \frac{c_1}{c_0}  \nabla u \cdot \nabla y_\alpha.
\]
Since
\[
\ds  \int_{\Omega} \nabla (u_\alpha -u) \cdot \nabla y_\alpha = -
\int_{\Omega}c_\alpha \partial_t (u_\alpha - u) (u_\alpha -u) = -
\frac{1}{2} \partial_t \int_{\Omega} c_\alpha(u_\alpha - u)^2,
\]
and
\[
\ds  \int_{\Omega} c_\alpha\partial_t^2 (u_\alpha -u) y_\alpha =
 - \frac{1}{2} \partial_t \int_{\Omega}
 |\nabla y_\alpha|^2,
\]
we obtain
\[
\ds
\partial_t \int_{\Omega}
 |\nabla y_\alpha|^2 +
\partial_t \int_{\Omega} c_\alpha(u_\alpha - u)^2 = -2 \alpha
\int_\Omega \frac{c_1}{c_0}  \nabla u \cdot \nabla y_\alpha \leq C
\alpha ||\nabla y_\alpha||_{L^\infty(0, T; L^2(\Omega))}.
\]
From the Gronwall Lemma it follows that \begin{equation} \label{u1}
|| u_\alpha - u ||_{L^\infty(0, T; L^2(\Omega))}  \leq C \alpha.
\end{equation} As a consequence, by using (\ref{hat-v}) one can see
that the function $\hat{u}_\alpha - \hat{u}$ solves the following
boundary value problem
\[\left\{
\begin{array}{l}
\Delta(\hat{u}_\alpha - \hat{u}) = O(\alpha) \quad {\rm in\;}
\Omega,\\ \nm (\hat{u}_\alpha - \hat{u}) |_{\partial \Omega} = 0.
\end{array}
\right.
\]
Integration by parts immediately gives, \begin{equation} \label{r4}
|| \grad(\hat{u}_\alpha - \hat{u})||_{L^2(\Omega)} = O(\alpha).
\end{equation} Taking into account that $\grad (u_\alpha - u) \in
L^{\infty}(0, T; L^2(\Omega))$, we find by using the above estimate
that
\begin{equation} \label{r42} || \grad(u_\alpha - u)||_{L^2(\Omega)}
= O(\alpha) {\rm \; a.e. \; } t \in (0, T). \end{equation}
 Under relation (\ref{tilde-v}), one can define the function $\tilde{y}_\alpha$ as solution of
\[
\left\{ \begin{array}{l} \tilde{y}_\alpha \in H^1_0(\Omega),\\ \nm
\Delta \tilde{y}_\alpha = c_\alpha\partial_t (\tilde{u}_\alpha -
u_\alpha) \quad {\rm in\;} \Omega.
\end{array}
\right.
\]
Integrating by parts immediately yields
\[
\displaystyle \int_{\Omega}c_\alpha\partial_{t}^{2}(\tilde{u}_\alpha
-
u_\alpha)\tilde{y}_\alpha=-1/2\partial_t\int_{\Omega}|\nabla\tilde{y}_\alpha|^2,
\]
and
\[
\displaystyle \int_{\Omega}\nabla(\tilde{u}_\alpha -
u_\alpha)\nabla\tilde{y}_\alpha=-1/2\partial_t\int_{\Omega}c_\alpha(\tilde{u}_\alpha
- u_\alpha)^2.
\]
To proceed with the proof of estimate (\ref{r5}), we firstly remark
that the function $\tilde{u}_\alpha$ given by (\ref{tilde-v}) is a
solution of
\[
\left\{
\begin{array}{l}
\ds (c_0\partial_t^2 - \Delta) \tilde{u}_\alpha  = i \alpha^d
\nabla\cdot (\eta c_1(x) e^{i \eta \cdot x})e^{-i|\eta|t} \in
L^2(\Omega) \quad {\rm in}\;
\Omega \times (0, T),\\
\nm \tilde{u}_\alpha  |_{t=0} = \varphi(x)  \quad {\rm in}\;
\Omega,\\
\nm \ds
 \partial_t \tilde{u}_\alpha
 |_{t=0} = \psi(x)  \quad {\rm in}\;
\Omega,\\ \nm \tilde{u}_\alpha  |_{\partial \Omega  \times (0, T)} =
e^{i \eta \cdot x - i |\eta| t}.
\end{array}
\right.
\]
Then, we deduce that $u_\alpha - \tilde{u}_\alpha$ solves the
following initial boundary value problem,
 \begin{equation}
\label{r2} \left\{
\begin{array}{l}
\ds (c_{\alpha}\partial_t^2 - \nabla\cdot~ \Delta ) (u_\alpha -
\tilde{u}_\alpha) = \alpha^d \nabla\cdot~ (c_1(x)\mbox{ grad}~
(\int_0^t e^{- i | \eta | s} v_{\eta}(x, t-s)\; ds) ) \quad {\rm
in}\; \Omega \times (0, T),\\ \nm
  (u_\alpha - \tilde{u}_\alpha)  |_{t=0} = 0 \quad {\rm in}\;
\Omega,\\
\nm \ds
 \partial_t  (u_\alpha - \tilde{u}_\alpha)
 |_{t=0} = 0  \quad {\rm in}\;
\Omega,\\
\nm
 (u_\alpha - \tilde{u}_\alpha)   |_{\partial \Omega  \times (0, T)}
= 0.
\end{array}
\right. \end{equation} Finally, we can use (\ref{r2}) to find by
integrating by parts that
\[
\ds
\partial_t \int_{\Omega}
 |\nabla \tilde{y}_\alpha|^2 +
\partial_t \int_{\Omega} c_\alpha(\tilde{u}_\alpha - u_\alpha)^2 = 2
\alpha^d \int_\Omega c_1  \grad (u - u_\alpha) \cdot \grad
\tilde{y}_\alpha
\]
which, from the Gronwall Lemma and by using (\ref{r42}), yields \[
|| \tilde{u}_\alpha - u_\alpha
 ||_{L^\infty(0, T; L^2(\Omega))}  \leq C^{\prime} \alpha^{d+1}.
\] This achieves the proof. $\square$\\

Now, we identify the function $ c(x)$ by using the difference
between local Dirichlet to Neumann maps and the function
$\theta_\eta$ as solution to the Volterra equation (\ref{eq4m}) or
equivalently the ODE (\ref{eq4p}),
as a function of $\eta$. Then, the following main result holds.\\

\begin{theorem} \label{th1}
Let $\eta \in \R^d$, $d=2,3$. Suppose that the smooth coefficient
$c(x)$ satisfies (\ref{coeff})-(\ref{c0-c1})-(\ref{ineq1}). Let
$u_\alpha$ be the unique solution in $ {\cal C}^0(0, T; H^1(\Omega))
\cap {\cal C}^1(0, T; L^2(\Omega))$ to the wave equation
(\ref{walpha}) with $ \varphi(x) = e^{i \eta \cdot x},$ $ \psi(x) =
- i | \eta | e^{i \eta \cdot x},$ and $ f(x, t) = e^{i \eta \cdot x
-  i | \eta | t}.$ Let $\tilde{f}$ be the function which satisfies
(\ref{f-tilde}). Suppose that $\Gamma$ and $T$ geometrically control
$\Omega$, then we have
\begin{equation} \label{eqa} \ds \int_0^T \int_\Gamma \Big(
\theta_\eta +
\partial_t \theta_\eta \partial_t \cdot\Big)(\Lambda_\alpha -\Lambda_0)( \tilde{f})(x,t)~d\sigma(x)dt = \alpha^{d-1}
|\eta|^2 \int_{\Omega^\prime} (c_\alpha-c_0)(x) e^{2 i \eta \cdot x}
\; dx + O(\alpha^{d+1}) \end{equation} \[=\alpha^{d} |\eta|^2
\int_{\Omega^\prime} c_1(x) e^{2 i \eta \cdot x} \; dx +
O(\alpha^{d+1}) ,\]where $\theta_\eta$ is the unique solution to the
ODE (\ref{eq4p}) with $g_\eta$ defined as the boundary control in
(\ref{weta}). The term $O(\alpha^{d+1})$ is independent of the
function $c_1$. It depends only on the bound $M$.
\end{theorem}

\proof Since the extension of $(\Lambda_\alpha -\Lambda_0)(
\tilde{f})(x,t)$ to $\partial\Omega\times(0,T)$ is
$\displaystyle(\frac{\partial u_\alpha}{\partial n} - \frac{\partial
u}{\partial n} )$, then by conditions $\partial_t \theta_\eta (T) =
0$ and $\displaystyle(\frac{\partial u_\alpha}{\partial n} -
\frac{\partial u}{\partial n} )|_{t=0}  = 0,$ we have
$(\Lambda_\alpha -\Lambda_0)( \tilde{f})(x,t)|_{t=0}=0.$ Therefore
the term $$ \ds \int_0^T \int_\Gamma
\partial_t \theta_\eta \partial_t  (\Lambda_\alpha -\Lambda_0)(
\tilde{f})(x,t)~d\sigma(x)dt$$ may be simplified as follows
\begin{equation} \label{d}  \ds \int_0^T \int_\Gamma
\partial_t \theta_\eta \partial_t  (\Lambda_\alpha -\Lambda_0)(
\tilde{f})(x,t)~d\sigma(x)dt= - \int_0^T \int_\Gamma
\partial^2_t \theta_\eta (\Lambda_\alpha -\Lambda_0)(
\tilde{f})(x,t)~d\sigma(x)dt. \end{equation}

On the other hand, we have
\[\begin{array}{l}
\ds \int_0^T \int_\Gamma
 \Bigr[
\theta_\eta (\Lambda_\alpha -\Lambda_0)( \tilde{f}) +
\partial_t \theta_\eta \partial_t (\Lambda_\alpha -\Lambda_0)(
\tilde{f}) \Bigr](x,t)d\sigma(x)dt =  \int_0^T \int_\Gamma
 \Bigr[
\theta_\eta (\Lambda_\alpha( \tilde{f})
-\tilde{\Lambda}_\alpha(\tilde{u}_\alpha|_{\Gamma\times(0,T)}))+\end{array}\]
\[\begin{array}{l}
\partial_t \theta_\eta \partial_t (\Lambda_\alpha( \tilde{f})
-\tilde{\Lambda}_\alpha(\tilde{u}_\alpha|_{\Gamma\times(0,T)}))
\Bigr](x,t)d\sigma(x)dt+\end{array}\]\[\begin{array}{l}
  \ds \int_0^T \int_{\Gamma} \Bigr[
\theta_\eta \alpha^d\int_0^t e^{-  i | \eta | s} \frac{\partial
v_{\eta}}{\partial n}(x, t-s)\; ds + \alpha^d\partial_t \theta_\eta
\partial_t \int_0^t e^{- i | \eta | s} \frac{\partial
v_{\eta}}{\partial n}(x, t-s)\; ds \Bigr]d\sigma(x)dt ;
\end{array}
\]
where
$\ds\tilde{\Lambda}_\alpha(\tilde{u}_\alpha|_{\Gamma\times(0,T)})=\Lambda_0(
\tilde{f})+\alpha^d\int_0^t e^{-  i | \eta | s} \Lambda_0(
v_{\eta}|_{\Gamma})(x, t-s)\; ds $.\\
Given that, $\theta_\eta$ satisfies the Volterra equation
(\ref{eq4p}) and
\[\begin{array}{l}
\ds \partial_t ( \int_0^t e^{-  i | \eta | s} \frac{\partial
v_{\eta}}{\partial n}(x, t-s)\; ds ) =
\partial_t (- e^{- i |\eta| t}  \int_0^t e^{ i | \eta | s} \frac{\partial
v_{\eta}}{\partial n}(x, s)\; ds ) \\ \nm \ds =  i |\eta| e^{- i
|\eta| t} \int_0^t e^{ i | \eta | s} \frac{\partial
v_{\eta}}{\partial n}(x, s)\; ds +  \frac{\partial
v_{\eta}}{\partial n}(x, t),
\end{array}
\]
we obtain by integrating by parts over $(0, T)$ that
\[\begin{array}{l}
  \ds \int_0^T \int_{\Gamma} \Bigr[
\theta_\eta \int_0^t e^{-  i | \eta | s} \frac{\partial
v_{\eta}}{\partial n}(x, t-s)\; ds + \partial_t \theta_\eta
\partial_t \int_0^t e^{- i | \eta | s} \frac{\partial
v_{\eta}}{\partial n}(x, t-s)\; ds \Bigr]d\sigma(x)dt  \\= \ds
\int_0^T \int_{\Gamma} \big( \frac{\partial v_{\eta}}{\partial n}(x,
t) (\partial_t \theta_\eta + \int_t^T \theta_\eta(s) e^{i |\eta|
(t-s)} \;ds) - i |\eta| (e^{- i |\eta| t} \partial_t \theta_\eta
(t)) \int_0^t e^{ i | \eta | s} \frac{\partial v_{\eta}}{\partial
n}(x, s)\; ds\big)\;d\sigma(x)dt\\ \nm = \ds  \int_0^T \int_{\Gamma}
\frac{\partial v_{\eta}}{\partial n}(x, t) (\partial_t \theta_\eta +
\int_t^T (\theta_\eta(s) - i |\eta| \partial_t \theta_\eta(s))
 e^{ i | \eta |(t- s)} \;ds) \;d\sigma(x)dt\\
\nm \ds
 = \int_0^T \int_{\Gamma} g_\eta(x, t)
\Lambda_0( v_{\eta}|_{\Gamma})(x, t)\; d\sigma(x)dt ,
\end{array}
\]
and so, from Proposition \ref{p4.1} we obtain
\[\begin{array}{l}
\ds \int_0^T \int_\Gamma
 \Bigr[
\theta_\eta (\Lambda_\alpha -\Lambda_0)( \tilde{f}) +
\partial_t \theta_\eta \partial_t (\Lambda_\alpha -\Lambda_0)(
\tilde{f}) \Bigr](x,t)d\sigma(x)dt  = \\ \nm \ds \alpha^d |\eta |^2
\int_{\Omega^\prime} c_1(x) e^{2 i \eta \cdot x} \; dx
\\
\nm \ds + \int_0^T \int_\Gamma
 \Bigr[
\theta_\eta (\Lambda_\alpha( \tilde{f})
-\tilde{\Lambda}_\alpha(\tilde{u}_\alpha|_{\Gamma\times(0,T)})) +
\partial_t \theta_\eta \partial_t (\Lambda_\alpha( \tilde{f})
-\tilde{\Lambda}_\alpha(\tilde{u}_\alpha|_{\Gamma\times(0,T)}))
\Bigr]d\sigma(x)dt + O(\alpha^{d+1}).
\end{array}
\]
Thus, to prove Theorem \ref{th1} it suffices then to show that
\[
 \int_0^T \int_\Gamma\Bigr[ \theta_\eta (\Lambda_\alpha( \tilde{f})
-\tilde{\Lambda}_\alpha(\tilde{u}_\alpha|_{\Gamma\times(0,T)})) +
\partial_t \theta_\eta \partial_t (\Lambda_\alpha( \tilde{f})
-\tilde{\Lambda}_\alpha(\tilde{u}_\alpha|_{\Gamma\times(0,T)}))
\Bigr]d\sigma(x)dt = O(\alpha^{d+1}).\]

From definition (\ref{hat-v}) we have
\[
\hat{u}_\alpha - \hat{\tilde{u}}_\alpha =  \int_0^T (u_\alpha -
\tilde{u}_\alpha) z(t)\; dt,\] which gives by system (\ref{r2}) that
\[
\Delta (\hat{u}_\alpha - \hat{\tilde{u}}_\alpha) =
\int_0^Tc_\alpha\partial_t^2 (u_\alpha - \tilde{u}_\alpha) z(t)\;
dt+\alpha^d\int_0^T\nabla\cdot(c_1(x)\mbox{
grad}~(\int_0^te^{-i|\eta|s}v_{\eta}(x,t-s)ds))z(t)~dt.
\]
Thus, by (\ref{tilde-v}) and (\ref{r2}) again, we see that the
function $\hat{u}_\alpha - \hat{\tilde{u}}_\alpha$  is solution of
 \begin{equation} \label{r3} \left\{
\begin{array}{l}
\ds - \Delta (\hat{u}_\alpha - \hat{\tilde{u}}_\alpha)  = -
\int_0^Tc_\alpha (u_\alpha - \tilde{u}_\alpha) z^{\prime
\prime}(t)\; dt + \nabla\cdot~ (c_1(x) \mbox{ grad}~
(\hat{\tilde{u}}_\alpha - \hat{u})) \quad {\rm in}\; \Omega,\\ \nm
 (\hat{u}_\alpha - \hat{\tilde{u}}_\alpha)   |_{\partial \Omega}
= 0.
\end{array}
\right. \end{equation} Taking into account estimate (\ref{r5}) given
by Lemma \ref{lemm1}, then by using standard elliptic regularity
(see e.g. \cite{E}) for the boundary value problem (\ref{r3}) we
find that
\[\ds || \frac{\partial} {\partial n} (\hat{u}_\alpha -
\hat{\tilde{u}}_{\alpha}) ||_{L^2(\Gamma)} = O(\alpha^{d+1}).\] The
fact that $\Lambda_\alpha( \tilde{f})
-\tilde{\Lambda}_\alpha\big(\tilde{u}_\alpha|_{\Gamma\times(0,T)}\big):=\displaystyle\frac{\partial}
{\partial n} (u_\alpha - \tilde{u}_{\alpha})\in
L^\infty(0,T;L^2(\Gamma))$, we deduce, as done in the proof of
Lemma~\ref{lemm1}, that
\[\ds || \Lambda_\alpha( \tilde{f})
-\tilde{\Lambda}_\alpha(\tilde{u}_\alpha|_{\Gamma\times(0,T)})||_{L^2(\Gamma)}
= O(\alpha^{d+1}),\]which implies that
\[
 \int_0^T \int_\Gamma\Bigr[ \theta_\eta (\Lambda_\alpha( \tilde{f})
-\tilde{\Lambda}_\alpha(\tilde{u}_\alpha|_{\Gamma\times(0,T)})) +
\partial_t \theta_\eta \partial_t (\Lambda_\alpha( \tilde{f})
-\tilde{\Lambda}_\alpha(\tilde{u}_\alpha|_{\Gamma\times(0,T)}))
\Bigr]d\sigma(x)dt = O(\alpha^{d+1}).\] This completes the proof of
our Theorem.
 $\square$\\


We are now in position to describe our identification procedure
which is based on Theorem \ref{th1}. Let us neglect the
asymptotically small remainder in the asymptotic formula
(\ref{eqa}). Then, it follows
\[ \ds c_\alpha(x)-c_0(x) \approx
\frac{2}{\alpha^{d-1}} \int_{\R^d} \frac{e^{-2 i \eta \cdot
x}}{|\eta|^2} \int_0^T \int_\Gamma \Big( \theta_\eta +
\partial_t \theta_\eta \partial_t \cdot\Big)(\Lambda_\alpha -\Lambda_0)( \tilde{f})(x,t)d\sigma(y)dtd\eta, x \in
\Omega.
\] The
method of reconstruction we propose here consists in sampling values
of
\[
\ds \frac{1}{|\eta|^2} \int_0^T \int_\Gamma\Big( \theta_\eta +
\partial_t \theta_\eta \partial_t \cdot\Big)(\Lambda_\alpha -
\Lambda_0)( \tilde{f})(x,t)d\sigma(x)dt
\]
at some discrete set of points $\eta$ and then calculating the
corresponding inverse Fourier transform.\\

In the following, it is not hard to prove the more convenient
approximation in terms of the values of local Dirichlet-to-Neumann
maps $\Lambda_\alpha$ and
$\Lambda_0$ at $\tilde{f}$.\\

\begin{corollary}
Let $\eta \in \R^d$ and let $\tilde{f}$ be defined by
(\ref{f-tilde}). Suppose that $\Gamma$ and $T$ geometrically control
$\Omega$, then we have the following approximation
\begin{equation} \label{emc} \ds c_\alpha(x)\approx c_0(x)  -
\frac{2}{\alpha^{d-1}} \int_{\R^d} \frac{e^{-2 i \eta \cdot
x}}{|\eta|^2} \int_0^T \int_\Gamma \Bigr[ e^{i |\eta| t}\partial_t (
e^{-i |\eta| t} g_\eta(y,t)) (\Lambda_\alpha -\Lambda_0)(
\tilde{f})(y, t) \Bigr]d\sigma(y)dtd\eta, x \in \Omega,
\end{equation} where the boundary control $g_\eta$
is defined by (\ref{weta}).
\end{corollary}

\proof The term $\ds \int_0^T \int_\Gamma
\partial_t \theta_\eta \partial_t  (\Lambda_\alpha -\Lambda_0)(
\tilde{f})(x,t)~d\sigma(x)dt$, given in Theorem \ref{th1}, has to be
interpreted as follows:
\[
\ds \int_0^T \int_\Gamma
\partial_t \theta_\eta \cdot \partial_t (\Lambda_\alpha -\Lambda_0)(
\tilde{f})(x,t)d\sigma(x)dt = -  \int_0^T \int_\Gamma
\partial^2_t \theta_\eta \cdot (\Lambda_\alpha -\Lambda_0)(
\tilde{f})(x,t)d\sigma(x)dt,
\]
because $\theta_\eta |_{t=T} = 0$ and $\displaystyle\partial_t
(\frac{\partial u_\alpha}{\partial n}  - \frac{\partial u}{\partial
n} ) |_{t =0} = 0$. In fact, in view of the ODE (\ref{eq4p}), the
term $\ds \int_0^T \int_\Gamma \Bigr[ \theta_\eta (\Lambda_\alpha
-\Lambda_0) +
\partial_t \theta_\eta \cdot\partial_t (\Lambda_\alpha -\Lambda_0) \Bigr]\tilde{f}(x, t)d\sigma(x)dt$ may be simplified after integration by parts
over $(0, T)$ and use of the fact that $\theta_\eta$ is the solution
to the ODE (\ref{eq4p}) to become
\[
\ds - \ds \int_0^T \int_\Gamma e^{i |\eta| t}
\partial_t ( e^{-i |\eta| t} g_\eta)\cdot (\Lambda_\alpha -\Lambda_0)(
\tilde{f})(x,t)d\sigma(x)dt.
\]
Then, the desired approximation is established. $\square$
\section{Conclusion}
The use of approximate formula (\ref{eqa}), including the difference
between the local Dirichlet to Neumann maps, represents a promising
approach to the dynamical identification and reconstruction of a
coefficient which is unknown in a bounded domain(but it is known
outside of this domain) for a class of hyperbolic PDE. We believe
that this method will yield a suitable approximation to the
dynamical identification of small conductivity ball (of the form $z
+ \alpha D$) in a homogeneous medium in $\R^d$ from the boundary
measurements. We will present convenable numerical implementations
for this investigation. This issue will be considered in a
forthcoming work.


\end{document}